\documentclass[12pt]{article}
\pdfoutput=1
\usepackage[utf8]{inputenc}
\usepackage{amsfonts,amsmath,amsthm,amssymb}
\usepackage{nicefrac}
\usepackage{enumitem}
\usepackage[colorlinks,allcolors=blue]{hyperref}
\usepackage{url}
\usepackage{breakcites}
\usepackage[ruled,vlined,noresetcount]{algorithm2e}
\SetKw{kwParams}{Parameters:}
\usepackage[margin=1in]{geometry}

\newtheorem{theorem}{Theorem}

\newtheorem{corollary}{Corollary}

\newtheorem{remark}{Remark}

\newtheorem{example}{Example}

\let\Pr\relax
\DeclareMathOperator\Pr{\mathbb{P}}

\DeclareMathOperator\E{\mathbb{E}}

\DeclareMathOperator\var{var}

\DeclareMathOperator\Geom{Geometric}

\newcommand{\KLr}[2]{ { \overline{D} \left({#1} \;\middle\Vert\; {#2}\right) } }

\newcommand{\btheta}{{\pmb{\theta}}}

\newcommand{\bphi}{{\pmb{\phi}}}

\newcommand{\calA}{{\mathcal A}}

\newcommand{\calF}{{\mathcal F}}

\newcommand{\Intpp}[1]{ { {\mathbb Z_{> 0}}^{#1} } }
\newcommand{\Real}[1]{ { {\mathbb R}^{#1} } }

\DeclareMathOperator*{\argmax}{arg\,max}

\title{\textbf{A Hoeffding Inequality for Finite State Markov Chains and its Applications to Markovian Bandits}}
\author{
Vrettos Moulos
\thanks{Supported in part by the NSF grant CCF-1816861.}
\\
University of California Berkeley \\
\href{mailto:vrettos@berkeley.edu}{vrettos@berkeley.edu}
}
\date{}

\begin{document}

\maketitle

\begin{abstract}
This paper develops a Hoeffding inequality for the partial sums $\sum_{k=1}^n f (X_k)$, where $\{X_k\}_{k \in \mathbb{Z}_{> 0}}$ is an irreducible
Markov chain on a finite state space $S$,
and $f : S \to [a, b]$ is a real-valued function.
Our bound is simple, general, since it only assumes irreducibility and finiteness of the state space, and powerful.
In order to demonstrate its usefulness we provide two applications in multi-armed bandit problems.
The first is about identifying an approximately best Markovian arm,
while the second is concerned with regret minimization in the context of Markovian bandits.
\end{abstract}

\section{Introduction}

Let $\{X_k\}_{k \in \Intpp{}}$ be a Markov chain on a finite state space $S$,
with initial distribution $q$, and irreducible transition probability matrix $P$, governed by the probability law $\Pr_q$.
Let $\pi$ be its stationary distribution, and $f : S \to [a, b]$ be a real-valued function on the state space. Then the strong law of large numbers for Markov chains asserts that,
\[
\frac{1}{n} \sum_{k=1}^n f (X_k) \stackrel{\Pr_q-a.s.}{\to} 
\E_\pi [f (X_1)],~\text{as}~ n \to \infty.
\]
Moreover, the central limit theorem for Markov chains provides a rate for this convergence,
\[
\sqrt{n}
\left(\frac{1}{n} \sum_{k=1}^n f (X_k) - \E_\pi [f (X_1)]\right)
\stackrel{d}{\to} N (0, \sigma^2),~\text{as}~ n \to \infty,
\]
where $\sigma^2 = \lim_{n \to \infty} \frac{1}{n} \var_q \left(\sum_{k=1}^n f (X_k)\right)$ is the limiting variance.

Those asymptotic results are insufficient in many applications which require
finite-sample estimates. One of the most central such application is the convergence of Markov chain Monte Carlo (MCMC)
approximation techniques~\cite{MCMC53},
where a finite-sample estimate is needed to bound the approximation error. 
Further applications include theoretical computer science and the approximation of the permanent~\cite{JSV01}, as well as statistical learning theory and multi-armed bandit problems~\cite{Moulos19-bandits-identification}.

Motivated by this discussion we provide a finite-sample Hoeffding inequality for finite Markov chains. In the special case that the random variables $\{X_k\}_{k \in \Intpp{}}$ are independent and identically distributed according to $\pi$, Hoeffding's classic inequality~\cite{Hoeffding63} states that,
\[
\Pr \left(
\left|\sum_{k=1}^n \left(f (X_k) - \E [f (X_1)]\right)\right| \ge  t
\right) \le
2 \exp\left\{- \frac{t^2}{2 \nu^2}\right\},
\]
where $\nu^2 = \frac{1}{4} n (b-a)^2$.
In Theorem~\ref{thm:Hoeffding} we develop a version of Hoeffding's inequality for finite state Markov chains. 
Our bound is very simple and easily computable, since it is based on martingale techniques and it only involves hitting times of Markov chains
which are very well studied for many types of Markov chains~\cite{Aldous-Fill-02}.
It is worth mentioning that our bound is based solely on irreducibility, 
and it does not make any extra assumptions like aperiodicity or reversibility which prior works require.

There is a rich literature on finite-sample bounds for Markov chains.
One of the earliest works~\cite{Davisson-Longo-Sgarro-81} uses counting and a generalization of the method of types, in order to derive a Chernoff bound for Markov chains which are irreducible and aperiodic. 
An alternative approach~\cite{WH17, Moulos-Ananth-19},
uses the theory of large deviations
to derive sharper Chernoff bounds.
When reversibility is assumed, the transition probability matrix is symmetric with respect to the space $L^2 (\pi)$, which enables the use of matrix perturbation theory.
This idea leads to Hoeffding inequalities that involve the spectral gap of the Markov chain and was initiated in~\cite{Gillman93}.
Refinements of this bound were given in a series of works~\cite{Dinwoodie95,Kahale97,Lezaud98,LP04,Miaso14}.
In~\cite{Paulin15,Rao19,FJS18-Hoeffding} a generalized spectral gap is introduced
in order to obtain bounds even for a certain class of irreversible Markov chains as long as they posses a strictly positive generalized spectral gap.
Information-theoretic ideas are used in~\cite{KLM06} in order to derive a Hoeffding inequality for Markov chains with general state spaces that satisfy Doeblin's minorization condition, which in the case of a finite state space can be written as,
\begin{equation}\label{eqn:Doeblin}
\exists m \in \Intpp{} ~ \exists y \in S ~ \forall x \in S: ~ P^m (x, y) > 0.
\end{equation}
Of course there are irreducible transition probability matrices $P$ for which~\eqref{eqn:Doeblin} fails, but if we further assume aperiodicity,
then~\eqref{eqn:Doeblin} is satisfied. Our approach uses Doob's martingale combined with Azuma's inequality,
and is probably closest related to the work in~\cite{Glynn-Ormoneir-02},
where a bound for Markov chains with general state spaces is established using Dynkin's martingale. But the result in~\cite{Glynn-Ormoneir-02} heavily relies on the Markov chains satisfying Doeblin's condition~\eqref{eqn:Doeblin}.
A regeneration approach for uniformly ergodic Markov chains, where one splits the Markov chain in i.i.d. blocks, and reduces the problem to the concentration of an i.i.d. process, can be found in~\cite{Douc-Moulines-Olsson-van-Handel-11}. 

Another line of research is related to the concentration of a function of $n$ random variables around its mean, under Markovian or other dependent structures.
This was pioneered by the works of Marton~\cite{Marton-96-a,Marton-96-b,Marton-98} who used the transportation method, and further developed using coupling ideas in~\cite{Samson-00,Marton-03,Chazottes-Collet-Kulske-Redig-07,Kontorovich-Ramanan-08,Paulin15}.

To illustrate the applicability of our bound we use it to study two Markovian multi-armed bandit problems. The stochastic multi-armed bandits problem is a prototypical statistical learning problem that exhibits an exploration-exploitation trade-off. One is given multiple options, referred to as arms, and each of them associated with a probability distribution. The emphasis is put on focusing as quickly as possible on the best available option, rather than estimating with high confidence the statistics of each option. The cornerstone of this field is the pioneering work of Lai and Robbins~\cite{Lai-Robbins-85}. Here we study two variants of the multi-armed bandits problem where the probability distributions of the arms form Markov chains. First we consider the task of identifying with some fixed confidence an approximately best arm, and we use our bound to analyze the median elimination algorithm, originally proposed in~\cite{EMM06} for the case of i.i.d. bandits. Then we turn into the problem of regret minimization for Markovian bandits, where we analyze the UCB algorithm that was introduced in~\cite{Auer-Cesa-Fischer-02} for i.i.d. bandits. For a thorough introduction to multi-armed bandits we refer the interested reader to the survey~\cite{BC12}, and the books~\cite{Lat-Szep-19,Slivkins-19}.
\section{A Hoeffding Inequality for Irreducible Finite State Markov Chains}

The central quantity that shows up in our Hoeffding inequality,
and makes it differ from the classical i.i.d. Hoeffding inequality,
is the maximum hitting time of a Markov chain with an irreducible transition probability matrix $P$.
This is defined as,
\[
\mathrm{HitT} (P) = \max_{x, y \in S}\: \E [T_y \mid X_1 = x],
\]
where $T_y = \inf \{n \ge 0 : X_{n+1} = y\}$ is the number of transitions taken in order to visit state $y$ for the first time.
$\mathrm{HitT} (P)$ is ensured to be finite due to irreducibility and the finiteness of the state space.

\begin{theorem}\label{thm:Hoeffding}
Let $\{X_k\}_{k \in \Intpp{}}$ be a Markov chain on a finite state space $S$, driven by an initial distribution $q$, and an irreducible transition probability matrix $P$. Let $f : S \to [a, b]$ be a real-valued function.
Then, for any $t > 0$,
\[
\Pr_q \left(
\left|\sum_{k=1}^n \left(f (X_k) - \E_q \left[f (X_k)\right]\right)\right|
\ge t \right)
\le
2 \exp\left\{-\frac{t^2}{2 \nu^2}\right\},
\]
where $\nu^2 = \frac{1}{4} n (b-a)^2 \mathrm{HitT} (P)^2$.
\begin{proof}
We define the sums $S_{l, m} = f (X_l) + \ldots + f (X_m)$, for $1 \le l \le m \le n$,
and the filtration $\calF_0 = \sigma (\emptyset), ~ \calF_k = \sigma (X_1, \ldots, X_k)$ for $k = 1, \ldots, n$.
Then $\left\{\E (S_{1, n} \mid \calF_k) - \E (S_{1, n} \mid \calF_0)\right\}_{k=0}^n$,
is a zero mean martingale with respect to $\{\calF_k\}_{k=0}^n$, and let
$\Delta_k = \E (S_{1, n} \mid \calF_k) - \E (S_{1, n} \mid \calF_{k-1})$,
for $k=1, \ldots, n$,
be the martingale differences. 

We first note the following bounds on the martingale differences,
\[
\min_{y \in S}\: \E (S_{1, n} \mid \calF_{k-1}, X_k = y) - \E (S_{1, n} \mid \calF_{k-1}) \le \Delta_k,
\]
and
\[
\Delta_k \le 
\max_{x \in S}\: \E (S_{1, n} \mid \calF_{k-1}, X_k = x) - \E (S_{1, n} \mid \calF_{k-1}).
\]
Therefore, in order to bound the variation of $\Delta_k$ it suffices to control,
\begin{align*}
& \max_{x \in S}\: \E (S_{1, n} \mid \calF_{k-1}, X_k = x) -
\min_{y \in S}\: \E (S_{1, n} \mid \calF_{k-1}, X_k = y) \\
&\qquad= \max_{x, y \in S}\:\left\{
\E [S_{k,n} \mid X_k = x] - \E [S_{k,n} \mid X_k = y]
\right\} \\
&\qquad= \max_{x, y \in S}\:\left\{
\E [S_{1,n-k+1} \mid X_1 = x] - \E [S_{1,n-k+1} \mid X_1 = y]
\right\},
\end{align*}
where in the first equality we used the Markov property, and in the second the time-homogeneity.

We now use a hitting time argument. 
Observe the following pointwise statements,
\[
S_{1,n-k+1} \le T_y b + S_{T_y+1, n-k+1}, \quad
S_{T_y+1, n-k+1} + T_y a \le S_{T_y+1, T_y + n-k+1},
\]
from which we deduce that,
\[
S_{1,n-k+1} \le T_y (b-a) + S_{T_y+1, T_y + n-k+1}.
\]
Taking $\E [ \cdot \mid X_1 = x]$-expectations, and using the strong Markov property we obtain,
\[
\E [ S_{1, n-k+1} \mid X_1 = x] \le 
(b-a) \E [T_y \mid X_1 = x] +
\E [S_{1, n-k+1} \mid X_1 = y].
\]
Therefore,
\[
\max_{x, y \in S}\:\left\{
\E [S_{1,n-k+1} \mid X_1 = x] - \E [S_{1,n-k+1} \mid X_1 = y]
\right\} \le 
(b-a) \mathrm{HitT} (P).
\]
With this in our possession we apply Hoeffding's lemma, see for instance Lemma 2.3 in~\cite{Luc-Lugosi-01},
in order to get,
\[
\E \left(e^{\theta \Delta_k} \mid \calF_{k-1}\right) \le \exp\left\{\frac{\theta^2 (b-a)^2 \mathrm{HitT} (P)^2}{8}\right\} = 
\exp\left\{\frac{\theta^2 \nu^2}{2 n}\right\},~\text{for all}~ \theta \in \Real{}.
\]
Using Markov's inequality, and successive conditioning we obtain that for $\theta > 0$,
\begin{align*}
\Pr \left(\sum_{k=1}^n \left(f (X_k) - \E_q \left[f (X_k)\right]\right) \ge t\right) 
&\le e^{-\theta t} \E \left[e^{\theta \left(\sum_{k=1}^n \Delta_k\right)}\right] \\
&= e^{-\theta t} \E \left[\E \left(e^{\theta \Delta_n} \middle| \calF_{n-1} \right) e^{\theta \left(\sum_{k=1}^{n-1} \Delta_k\right)}\right] \\
&\le \exp{\left\{-\theta t + \frac{\theta^2 \nu^2}{2 n}\right\}}
\E \left[e^{\theta \left(\sum_{k=1}^{n-1} \Delta_k\right)}\right] \\
&\le \ldots \le \exp\left\{-\theta t + \frac{\theta^2 \nu^2}{2}\right\}.
\end{align*}
Plugging in $\theta = t/\nu^2$, we see that,
\[
\Pr \left(\sum_{k=1}^n \left(f (X_k) - \E_q \left[f (X_k)\right]\right) \ge t\right) \le \exp\left\{-\frac{t^2}{2 \nu^2}\right\}.
\]
The conclusion follows by combining the inequality above for $f$ and $-f$.
\end{proof}
\end{theorem}

\begin{example}
Consider a two-state Markov chain with $S = \{0, 1\}$ and
$P (0, 1) = p, ~ P (1, 0) = r$, with $p, r \in (0, 1]$. Then,
\[
\mathrm{HitT} (P) = \max \{\E [\Geom (p)], \E [\Geom (r)]\} = 1/\min\{p, r\},
\]
and Theorem~\ref{thm:Hoeffding} takes the form,
\[
\Pr_q \left(
\left|\sum_{k=1}^n \left(f (X_k) - \E_q \left[f (X_k)\right]\right)\right|
\ge t \right)
\le
2 \exp\left\{-\frac{2 \min \{p^2, r^2\} t^2}{n (b-a)^2}\right\}.
\]
\end{example}

\begin{example}
Consider the random walk on the $m$-cycle with state space $S = \{0, 1, \ldots, m-1\}$, and transition probability matrix $P (x, y) = (1 \{y \equiv x+1 \pmod{m}\} + 1 \{y \equiv x-1 \pmod{m}\})/2$. If $m$ is odd, then the Markov chain is aperiodic, while if $m$ is even, then the Markov chain has period $2$. Then,
\[
\mathrm{HitT} (P) =
\max_{y \in S}\: \E [T_y \mid X_1 = 0] =
\max_{y \in S}\: y (m - y) = 
\lfloor m^2/4\rfloor,
\]
and Theorem~\ref{thm:Hoeffding} takes the form,
\[
\Pr_q \left(
\left|\sum_{k=1}^n \left(f (X_k) - \E_q \left[f (X_k)\right]\right)\right|
\ge t \right)
\le
2 \exp\left\{-\frac{2 t^2}{n (b-a)^2 \lfloor m^2/4\rfloor^2}\right\}.
\]
\end{example}

\begin{remark}
Observe that the technique used to establish Theorem~\ref{thm:Hoeffding}
is limited to Markov chains with a finite state space $S$.
Indeed, if $\{X_k\}_{k \in \Intpp{}}$ is a Markov chain on a countably infinite state space $S$ with an irreducible and positive recurrent transition probability matrix $P$ and a stationary distribution $\pi$, then
we claim that,
\[
\frac{1}{\pi (y)} \le 1 + 
\sup_{x \in S}\: \E [T_y \mid X_1 = x],
~\text{for all}~ y \in S,
\]
from which it follows that
$\sup_{x, y \in S} \E [T_y \mid X_1 = x] = \infty$,
due to the fact that $\sum_{y \in S} \pi (y) = 1$ and $S$ is countably infinite. The aforementioned inequality can be established as follows.
\begin{align*}
\frac{1}{\pi (y)}
&= \E [\inf \{n \ge 1 : X_{n+1} = y \} \mid X_1 = y] \\
&= \sum_{x \in S} \E [\inf \{n \ge 1 : X_{n+1} = y \} \mid X_2 = x] P (y, x) \\
&\le \sup_{x \in S}\: \E [\inf \{n \ge 1 : X_{n+1} = y\} \mid X_2 = x] \\
&= 1 + \sup_{x \in S}\: \E [T_y \mid X_1 = x].
\end{align*}
\end{remark}

Moreover, through Theorem~\ref{thm:Hoeffding} we can obtain a concentration inequality for sums of a function evaluated on the transitions of a Markov chain. In particular, let 
\[
S^{(2)} = \{(x, y) \in S \times S : P (x, y) > 0\}.
\]
On the state space $S^{(2)}$ define the transition probability matrix,
\[
P^{(2)} \left((x, y), (z, w)\right) = I \{y = z\} P (y, w),
~\text{for}~ (x, y), (z, w) \in S^{(2)}.
\]
It is straightforward to verify that the fact that $P$ is irreducible, implies that $P^{(2)}$ is irreducible as well. This readily gives the following theorem.
\begin{theorem}\label{thm:pair-Hoeffding}
Let $\{X_k\}_{k \in \Intpp{}}$ be a Markov chain on a finite state space $S$, driven by an initial distribution $q$, and an irreducible transition probability matrix $P$. Let $f^{(2)} : S^{(2)} \to [a, b]$ be a real-valued function evaluated on the transitions of the Markov chain.
Then, for any $t > 0$,
\[
\Pr_q \left(
\left|\sum_{k=1}^n \left(f^{(2)} (X_k, X_{k+1}) - \E_q \left[f^{(2)} (X_k, X_{k+1})\right]\right)\right|
\ge t \right)
\le
2 \exp\left\{-\frac{t^2}{2 \nu^2}\right\},
\]
where $\nu^2 = \frac{1}{4} n (b-a)^2 \mathrm{HitT} \left(P^{(2)}\right)^2$.
\end{theorem}

\begin{corollary}\label{cor:lln}
When the Markov chain is initialized with its stationary distribution, $\pi$, Theorem~\ref{thm:Hoeffding} and Theorem~\ref{thm:pair-Hoeffding} give the following nonasymptotic versions of the weak law of large numbers for irreducible Markov chains. For any $\epsilon > 0$,
\[
\Pr_\pi \left(
\left|\frac{1}{n} \sum_{k=1}^n f (X_k) - \E_\pi \left[f (X_1)\right]\right|
\ge \epsilon \right)
\le
2 \exp\left\{-\frac{2 n \epsilon^2}{(b-a)^2 \mathrm{HitT} (P)^2}\right\},
\]
and,
\[
\Pr_\pi \left(
\left|\frac{1}{n} \sum_{k=1}^n f^{(2)} (X_k, X_{k+1}) - \E_\pi \left[f^{(2)} (X_1, X_2)\right]\right|
\ge \epsilon \right)
\le
2 \exp\left\{-\frac{2 n \epsilon^2}{(b-a)^2 \mathrm{HitT} \left(P^{(2)}\right)^2}\right\}.
\]
\end{corollary}
\section{Markovian Multi-Armed Bandits}

\subsection{Setup}

There are $K \ge 2$ arms, and each arm $a \in [K] = \{1, \ldots, K\}$
is associated with a parameter $\theta_a \in \Real{}$ which uniquely
encodes\footnote{$\Real{}$ and the set of $|S| \times |S|$ irreducible transition probability matrices have the same cardinality, and hence there is a bijection between them.} an irreducible transition probability matrix $P_{\theta_a}$.
We will denote the overall parameter configuration of all $K$ arms with $\btheta = (\theta_1, \ldots, \theta_K) \in \Real{K}$.
Arm $a$ evolves according to the stationary Markov chain, $\{X_n^a\}_{n \in \Intpp{}}$, driven by the irreducible transition probability matrix $P_{\theta_a}$ which has a unique stationary distribution $\pi_{\theta_a}$,
so that $X_1^a \sim \pi_{\theta_a}$.
There is a common reward function $f : S \to [c, d]$
which generates the reward process $\{Y_n^a\}_{n \in \Intpp{}} = \{f (X_n^a)\}_{n \in \Intpp{}}$.
The reward process, in general, is not going to be a Markov chain, unless $f$ is injective,
and it will have more complicated dependencies than the underlying Markov chain. Each time that we select arm $a$, this arm evolves by one transition and we observe the corresponding sample from the reward process $\{Y_n^a\}_{n \in \Intpp{}}$, while all the other arms stay rested.

The stationary reward of arm $a$ is $\mu (\theta_a) = \sum_{x \in S} f (x) \pi_{\theta_a} (x)$. 
Let $\mu^* (\btheta) = \max_{a \in [K]} \mu (\theta_a)$ be the maximum stationary mean, and for simplicity assume that there exists a unique arm,
$a^* (\btheta)$, attaining this maximum stationary mean, i.e. $\{a^* (\btheta)\} = \argmax_{a \in [K]} \mu (\theta_a)$.
In the following sections we will consider two objectives:
identifying an $\epsilon$ best arm with some fixed confidence level $\delta$ using as few samples as possible, and minimizing the expected regret given some fixed time horizon $T$.

\subsection{Approximate Best Arm Identification}

In the approximate best arm identification problem,
we are given an approximation accuracy $\epsilon > 0$,
and a confidence level $\delta \in (0, 1)$.
Our goal is to come up with an adaptive algorithm $\calA$ which 
collects a total of $N$ samples, and returns an arm $\hat{a}$
that is within $\epsilon$ from the best arm, $a^* (\btheta)$, with probability at least $1-\delta$, i.e.
\[
\Pr_\btheta^\calA (\mu^* (\btheta) \ge \mu (\theta_{\hat{a}}) + \epsilon) \le \delta.
\]
Such an algorithm is called $(\epsilon, \delta)$-PAC
(probably approximately correct).

In~\cite{Mannor-Tsitsiklis-04} a lower bound for the sample complexity of any $(\epsilon, \delta)$-PAC algorithm is derived. 
The lower bound states that no matter the $(\epsilon, \delta)$-PAC algorithm $\calA$, there exists an instance $\btheta$ such that the sample complexity is at least,
\[
\E_\btheta^\calA [N] =
\Omega \left(\frac{K}{\epsilon^2} \log \frac{1}{\delta}\right).
\]

A matching upper bound is provided for i.i.d. bandits in~\cite{EMM06}
in the form of the median elimination algorithm.
We demonstrate the usefulness of our Hoeffding inequality,
by providing an analysis of the median elimination algorithm 
in the more general setting of Markovian bandits.
\begin{algorithm}[ht]
\SetAlgoLined
\kwParams{
number of arms $K \ge 2$,
approximation accuracy $\epsilon > 0$,
confidence level $\delta \in (0, 1)$,
parameter $\beta$
}\;
$r = 1,~ A_r = [K],~ \epsilon_r = \epsilon/4,~ \delta_r = \delta/2$\;
\While{$|A_r| \ge 2$}{
$N_r = \left\lceil\frac{4 \beta}{\epsilon_r^2} \log \frac{3}{\delta_r}\right\rceil$\;
Sample each arm in $A_r$ for $N_r$ times\;
For $a \in A_r$ calculate $\Bar{Y}_a [r] = \frac{1}{N_r} \sum_{n=1}^{N_r} Y_n^a$\;
$m_r = \pmb{\mathrm{median}} \left((\Bar{Y}_a [r])_{a \in A_r}\right)$\;
Pick $A_{r+1}$ such that: 
\begin{itemize}
    \item $A_{r+1} \subseteq \{a \in A_r : \Bar{Y}_a [r] \ge m_r\}$\;
    \item $|A_{r+1}| = \lfloor |A_r|/2 \rfloor$\;
\end{itemize}
$\epsilon_{r+1} = 3 \epsilon_r / 4,~ \delta_{r+1} = \delta_r / 2,~ r = r + 1$\;
}
\Return{$\hat{a}$, where $A_r = \{\hat{a}\}$}\;
\caption{The $\beta$-Median-Elimination algorithm.}\label{alg:median}
\end{algorithm}

\begin{theorem}
If $\beta \ge \frac{1}{2} (d-c)^2 \max_{a \in [K]} \mathrm{HitT} (P_{\theta_a})^2$ then,
the $\beta$-Median-Elimination algorithm is $(\epsilon, \delta)$-PAC, and its sample complexity is upper bounded by $O \left(\frac{K}{\epsilon^2} \log \frac{1}{\delta}\right)$.
\begin{proof}
The total number of sampling rounds is at most $\lceil \log_2 K \rceil$,
and we can set them equal to $\lceil \log_2 K \rceil$ by setting $A_r = \{\hat{a}\}$, for $r \ge R_0$, where $A_{R_0} = \{\hat{a}\}$.
Fix $r \in \{ 1, \ldots, \lceil \log_2 K \rceil \}$. We claim that,
\begin{equation}\label{eqn:per-round}
\Pr_\btheta^{\beta-\mathrm{ME}} \left(\max_{a \in A_r} \mu (\theta_a) \ge \max_{a \in A_{r+1}} \mu (\theta_a) + \epsilon_r\right)
\le \delta_r.
\end{equation}
We condition on the value of $A_r$. If $|A_r| = 1$, then the claim is trivially true, so we only consider the case $|A_r| \ge 2$.
Let $\mu_r^* = \max_{a \in A_r} \mu (\theta_a)$, and $a_r^* \in \argmax_{a \in A_r : \mu (\theta_a) = \mu_r^*} \Bar{Y}_a [r]$.
We consider the following set of bad arms,
\[
B_r = \{
b \in A_r : \Bar{Y}_b [r] \ge \Bar{Y}_{a_r^*} [r],
~ \mu_r^* \ge \mu (\theta_b) + \epsilon_r
\},
\]
and observe that,
\begin{equation}\label{eqn:round-to-bad}
\Pr_\btheta^{\beta-\mathrm{ME}} \left(\mu_r^* \ge \mu_{r+1}^* + \epsilon_r\right)
\le \Pr_\btheta^{\beta-\mathrm{ME}} (|B_r| \ge  |A_r|/2).
\end{equation}
In order to upper bound the latter fix $b \in A_r$ and write,
\begin{align*}
& \Pr_\btheta^{\beta-\mathrm{ME}} \left(
\Bar{Y}_b [r] \ge \Bar{Y}_{a_r^*} [r],
\mu_r^* \ge \mu (\theta_b) + \epsilon_r
\middle|
\Bar{Y}_{a_r^*} [r] > \mu_r^* - \epsilon_r/2\right) \\
& \qquad \le \Pr_{\theta_b} (\Bar{Y}_b [r] \ge \mu (\theta_b) + \epsilon_r/2)
\le \delta_r/3,
\end{align*}
where in the last inequality we used~\autoref{cor:lln}.
Now via Markov's inequality this yields,
\begin{equation}\label{eqn:bound-bad}
\Pr_\btheta^{\beta-\mathrm{ME}} \left(|B_r| \ge  |A_r|/2
\middle|
\Bar{Y}_{a_r^*} [r] > \mu_r^* - \epsilon_r/2\right) \\
\le 2 \delta_r/3.
\end{equation}
Furthermore,~\autoref{cor:lln} gives that for any $a \in A_r$,
\begin{equation}\label{eqn:apply-Hoeffding}
\Pr_{\theta_a} (\Bar{Y}_a [r] \le \mu (\theta_a) - \epsilon_r/2)
\le \delta_r/3.
\end{equation}
We obtain~\eqref{eqn:per-round} by using~\eqref{eqn:bound-bad} and~\eqref{eqn:apply-Hoeffding} in~\eqref{eqn:round-to-bad}.

With~\eqref{eqn:per-round} in our possession,
the fact that median elimination is $(\epsilon, \delta)$-PAC
follows through a union bound,
\begin{align*}
\Pr_\btheta^{\beta-\mathrm{ME}} (\mu^* (\btheta) \ge \mu (\theta_{\hat{a}}) + \epsilon)
&\le \Pr_\btheta^{\beta-\mathrm{ME}} \left( 
\bigcup_{r=1}^{\lceil \log_2 K \rceil} \left\{
\mu_r^* \ge \mu_{r+1}^* + \epsilon_r
\right\}\right) \\
&\le \sum_{r=1}^\infty \delta_r \le \delta.
\end{align*}

Regarding the sample complexity, we have that the total number of samples is at most,
\begin{align*}
K \sum_{r=1}^{\lceil \log_2 K \rceil} N_r/2^{r-1}
&\le 2 K + \frac{64 \beta K}{\epsilon^2}
\sum_{r=1}^\infty \left(\frac{8}{9}\right)^{r-1} \log \frac{ 2^r 3}{\delta} \\
&= O \left(\frac{K}{\epsilon^2} \log \frac{1}{\delta} \right). 
\end{align*}
\end{proof}
\end{theorem}
\subsection{Regret Minimization}

Our device to solve the regret minimization problem is
an \emph{adaptive allocation rule},
$\bphi = \{\phi_t\}_{t \in \Intpp{}}$,
which is a sequence of random variables
where $\phi_t \in [K]$ is the arm that we select at time $t$.
Let $N_a (t) = \sum_{s=1}^t I_{\{\phi_s = a\}}$, be the number of times
we selected arm $a$ up to time $t$.
Our decision, $\phi_t$, at time $t$ is based on the information that we have accumulated so far.
More precisely, the event $\{\phi_t = a\}$ is measurable with respect to
the $\sigma$-field generated by the past decisions $\phi_1, \ldots, \phi_{t-1}$, and the past observations $\{X_n^1\}_{n=1}^{N_1 (t-1)}, \ldots, \{X_n^K\}_{n=1}^{N_K (t-1)}$.

Given a time horizon $T$, and
a parameter configuration $\btheta$,
the expected regret incurred when the adaptive allocation rule $\bphi$ is used,
is defined as,
\[
R_\btheta^\bphi (T) =
\sum_{b \not\in a^* (\btheta)}
\E_\btheta^\bphi [N_b (T)] \Delta_b (\btheta),
\]
where $\Delta_b (\btheta) = \mu^* (\btheta) - \mu (\theta_b)$.
Our goal is to come up with an adaptive allocation rule that makes the expected regret as small as possible.

There is a known asymptotic lower bound on how much we can minimize the expected regret.
Any adaptive allocation rule that is uniformly good across all parameter configurations should satisfy the following instance specific, asymptotic regret lower bound (see~\cite{Ananth-Varaiya-Walrand-II-87} for details),
\[
\sum_{b \neq a^* (\btheta)}
\frac{\Delta_b (\btheta)}{\KLr{\theta_b}
{\theta_{a^* (\btheta)}}} \le 
\liminf_{T \to \infty} \frac{R_\btheta^\bphi (T)}{\log T},
\]
where $\KLr{\theta}{\lambda}$ is the Kullback-Leibler divergence rate 
between the Markov chains with transition probability matrices $P_\theta$ and $P_\lambda$, given by,
\[
\KLr{\theta}{\lambda} = 
\sum_{x, y \in S} \log \frac{P_\theta (x, y)}{P_\lambda (x, y)} \pi_\theta (x) P_\theta (x, y).
\]

Here we utilize our Theorem~\ref{thm:Hoeffding} to provide a finite-time analysis of the $\beta$-UCB adaptive allocation rule for Markovian bandits, which is order optimal. The $\beta$-UCB adaptive allocation rule, is a simple and computationally efficient index policy based on upper confidence bounds which was initially proposed in~\cite{Auer-Cesa-Fischer-02} for i.i.d. bandits. It has already been studied in the context of Markovian bandits in~\cite{Tekin-Liu-10}, but in a more restrictive setting under the further assumptions of aperiodicity and reversibility due to the use of the bounds from~\cite{Gillman93,Lezaud98}.
For adaptive allocation rules that asymptotically match the lower bound we refer the interested reader to~\cite{Ananth-Varaiya-Walrand-II-87,Moulos20-bandits-regret}.
\begin{algorithm}[ht]
\SetAlgoLined
\kwParams{
number of arms $K \ge 2$,
time horizon $T \ge K$,
parameter $\beta$
}\;
Pull each arm in $[K]$ once\;
\For{$t = K$ to $T-1$,}{
\[
\phi_{t+1} \in \argmax_{a \in [K]} \left\{ 
\Bar{Y}_a (t) +
\sqrt{\frac{2 \beta \log t}{N_a (t)}}
\right\}\;
\]
}
\caption{The $\beta$-UCB adaptive allocation rule.}\label{alg:UCB}
\end{algorithm}

\begin{theorem}
If $\beta > \frac{1}{2} (d-c)^2 \max_{a \in [K]} \mathrm{HitT} (P_{\theta_a})^2$ then,
\[
R_\btheta^{\bphi_{\beta-\mathrm{UCB}}} (T)
\le 
8 \beta \left(\sum_{b \neq a^* (\btheta)}
\frac{1}{\Delta_b (\btheta)}\right) \log T 
+
\frac{\gamma}{\gamma - 2} \sum_{b \neq a^* (\btheta)} \Delta_b (\btheta),
\]
where $\gamma = \frac{4 \beta}{(d-c)^2 \max_{a \in [K]} \mathrm{HitT} (P_{\theta_a})^2} > 2$.
\begin{proof}
Fix $b \neq a^* (\btheta)$, and observe that,
\[
N_b (T) \le 
1 + \frac{8 \beta}{\Delta_b (\btheta)^2} \log T +
\sum_{t=2}^{T-1} I_{\left\{\phi_{t+1} = b, ~ N_b (t) \ge \frac{8 \beta}{\Delta_b (\btheta)^2} \log T\right\}}.
\]
On the event $\left\{\phi_{t+1} = b, ~ N_b (t) \ge \frac{8 \beta}{\Delta_b (\btheta)^2} \log T\right\}$, we have that, either 
$\Bar{Y}_b (t) \ge \mu (\theta_b) + \sqrt{\frac{2 \beta \log t}{N_b (t)}}$,
or $\Bar{Y}_{a^* (\btheta)} (t) \le \mu^* (\btheta) - \sqrt{\frac{2 \beta \log t}{N_{a^* (\btheta)} (t)}}$, since otherwise the $\beta$-UCB index of $a^* (\btheta)$ is larger than the $\beta$-UCB index of $b$ which contradicts the assumption that $\phi_{t+1} = b$.

In addition, using~\autoref{cor:lln}, we obtain,
\begin{align*}
& \Pr_\btheta^{\bphi_{\beta - \mathrm{UCB}}} \left(
\Bar{Y}_b (t) \ge \mu (\theta_b) + \sqrt{\frac{2 \beta \log t}{N_b (t)}}
\right) \\
& \quad=
\sum_{n=1}^t \Pr_\btheta^{\bphi_{\beta - \mathrm{UCB}}} \left(
\Bar{Y}_b (t) \ge \mu (\theta_b) + \sqrt{\frac{2 \beta \log t}{N_b (t)}},
N_b (t) = n
\right) \\
& \quad\le 
\sum_{n=1}^t \Pr_{\theta_b} \left(
\frac{1}{n} \sum_{k=1}^n Y_k^b \ge \mu (\theta_b) + \sqrt{\frac{2 \beta \log t}{n}}
\right) \\
& \quad\le
\sum_{n=1}^t \frac{1}{t^\gamma} = \frac{1}{t^{\gamma-1}}.
\end{align*}
Similarly we can see that,
\[
\Pr_\btheta^{\bphi_{\beta - \mathrm{UCB}}} \left(
\Bar{Y}_{a^* (\btheta)} (t) \le \mu^* (\btheta) - \sqrt{\frac{2 \beta \log t}{N_{a^* (\btheta)} (t)}}
\right)
\le \frac{1}{t^{\gamma-1}}.
\]
The conclusion now follows by putting everything together and using the integral estimate,
\[
\sum_{t=2}^{T-1} \frac{1}{t^{\gamma - 1}} \le
\int_1^\infty \frac{1}{t^{\gamma - 1}} d t =
\frac{1}{\gamma-2}.
\]

\end{proof}
\end{theorem}

\section*{Acknowledgements}
This research was supported in part by the NSF grant CCF-1816861.

\bibliographystyle{apalike}
\bibliography{references}

\end{document}